\documentclass[10pt,a4paper]{amsart}

\usepackage[latin1]{inputenc}
\usepackage[T1]{fontenc}
\usepackage{eucal,url,amssymb,bm,stmaryrd,booktabs,verbatim,enumerate}
\usepackage[pagebackref]{hyperref}

\theoremstyle{plain}
\newtheorem{theorem}{Theorem}

\newtheorem{lemma}[theorem]{Lemma}
\newtheorem{corollary}[theorem]{Corollary}

\theoremstyle{definition}

\theoremstyle{remark}
\newtheorem{remark}[theorem]{Remark}

\newtheorem*{acknowledgements}{Acknowledgements}

\numberwithin{equation}{section}

\newcommand{\Lie}[1]{\operatorname{\textsl{#1}}}
\newcommand{\lie}[1]{\operatorname{\mathfrak{#1}}}

\newcommand{\SP}{\Lie{Sp}}

\newcommand{\SU}{\Lie{SU}}
\newcommand{\su}{\lie{su}}

\newcommand{\om}[1]{{\omega_{#1}}}

\newcommand{\rep}[2]{{\Lambda^{#1}_{#2}}}

\newcommand{\Id}{{\Lie{Id}}}

\newcommand{\LC}{\nabla^g}
\newcommand{\RC}{R^g}

\DeclareMathOperator{\vol}{\Lie{vol}}

\newcommand{\Hodge}{\mathord{\mkern1mu *}}
\newcommand{\hook}{\mathbin{\lrcorner}}
\newcommand{\bw}{\mathbin{\wedge}}

\newcommand{\abs}[1]{\left\lvert #1\right\rvert}
\newcommand{\norm}[1]{\left\lVert #1\right\rVert}

\newcommand{\real}[1]{\left\llbracket #1 \right\rrbracket}

\newcommand{\set}[2]{\left\{ #1\colon#2 \right\}}

\begin{document}

\title[Conformal equivalences in dimension 6 and 7]{Conformal
  equivalence between certain geometries in dimension $\mathbf 6$ and
  $\mathbf 7$}

\date{\today}

\author{Richard Cleyton}
\address[Cleyton]{Humboldt-Universit\"at zu Berlin\\
Institut f\"ur Mathematik\\
Unter den Linden 6\\
D-10099 Berlin\\
Germany}
\email{cleyton@mathematik.hu-berlin.de}

\author{Stefan Ivanov}
\address[Ivanov]{University of Sofia "St. Kl. Ohridski"\\
Faculty of Mathematics and Informatics,\\
Blvd. James Bourchier 5,\\
1164 Sofia, Bulgaria}
\email{ivanovsp@fmi.uni-sofia.bg}

\begin{abstract}
  For $G_2$-manifolds the Fern\'andez-Gray class $\mathcal
  X_1+\mathcal X_4$ is shown to consist of the union of the class
  $\mathcal X_4$ of $G_2$-manifolds locally conformal to parallel
  $G_2$-structures and that of conformal transformations of nearly
  parallel or weak holonomy $G_2$-manifolds of type $\mathcal X_1$.
  The analogous conclusion is obtained for Gray-Hervella class
  $\mathcal W_1+\mathcal W_4$ of real $6$-dimensional almost Hermitian
  manifolds: this sort of geometry consists of locally conformally
  K\"ahler manifolds of class $\mathcal W_4$ and conformal
  transformations of nearly K\"ahler manifolds in class $\mathcal
  W_1$.  A corollary of this is that a compact $\SU(3)$-space in class
  $\mathcal W_1+\mathcal W_4$ or $G_2$-space of the kind $\mathcal
  X_1+\mathcal X_4$ has constant scalar curvature if only if it is
  either a standard sphere or a nearly parallel $G_2$ or nearly
  K\"ahler manifold, respectively.  The properties of the Riemannian
  curvature of the spaces under consideration are also explored.
\end{abstract}
\maketitle
\section{Introduction}

Reductions of the bundle of orthonormal frames over a Riemannian
manifold to a principal $G$-bundle may be classified by the
$G$-invariant components of the intrinsic torsion.

This idea was originally due to Gray and collaborators~\cite{MR696037,
  MR581924} for the special instances of $G_2$-manifolds and almost
Hermitian manifolds.  It has been further refined and explored by, for
instance, Bryant~\cite{math.DG/0305124}, Farinola, Falcitelli \&
Salamon~\cite{MR1299398}, Mart\'in Cabrera~\cite{MR2162420,MR1386249},
Mart\'in Cabrera, Monar \& Swann~\cite{MR1373070}, Chiossi \&
Salamon~\cite{MR1922042}.

For $G_2$- and almost Hermitian structures alike, the intrinsic
torsion has 4 irreducible components. There are thus potentially 16
torsion classes for these two kinds of geometries.

In~\cite{MR1373070}, Mart\'in Cabrera, Monar \& Swann showed that
apart from one instance, $\mathcal X_1+\mathcal X_2$ in our notation,
every single class of $G_2$-structures may be realized on a compact
homogeneous space.  For the one exception an easy calculation shows
that any $G_2$ structure with torsion $X_1+X_2$ must have either
$X_1=0$ or $X_2=0$.

In section~\ref{sec:g_2-case} we show that something similar holds for
the class $\mathcal X_1+\mathcal X_4$.  Namely that the latter
essentially is generated by the classes $\mathcal X_1$ and $\mathcal
X_4$.  It is well known that $G_2$-structures in this class are
locally conformally equivalent to nearly parallel ones.  We will show
that this equivalence is only really local when the $G_2$-structure
lies in the subclass $\mathcal X_4$ of locally conformally parallel
structures. The structure of compact locally conformally parallel
$G_2$-manifolds has been recently described in~\cite{IPP,Verb}. In
contrast to this, we will show that if the $X_1$ component is non-zero
at some point, it is non-zero everywhere. This is the key point in
proving that a global conformal change exists.

Differeently from the $G_2$ case, it was only recently pointed out by
Butruille~\cite{math.DG/0503150} that a $6$-dimensional almost
Hermitian manifold in the Gray-Hervella class $\mathcal W_1+\mathcal
W_4$ is locally conformal to a nearly K\"ahler manifold.  In
section~\ref{sec:4}, we present a different and simpler proof of this
fact for completeness.  Based on this the analogous statements to
those given for $G_2$-structures are shown to hold for $6$-dimensional
almost Hermitian geometry, too. In particular, any almost Hermitian
$6$-manifold in the class $(\mathcal W_1+\mathcal
W_4)\setminus\mathcal W_4$ has trivial canonical bundle.  The
geometries $\mathcal W_4$ for $\SU(3)$- and $\mathcal X_4$ for
$G_2$-manifolds are both special instances of $G$-structures with
vectorial torsion.  This notion was studied in~\cite{math.DG/0509147}.
The almost Hermitian manifolds and $G_2$-structures studied in this
paper all fit in the wider framework of $G$-structures with three-form
torsion, see~\cite{FI,MR2047649}. 

The aim of this note is to establish
\begin{theorem}\label{main1}
  Let $(M,g,\phi)$ be a compact 7-dimensional manifold locally
  conformally equivalent to a nearly parallel $G_2$-manifold.  Then
  $(M,g,\phi)$ has constant scalar curvature if and only if $(M,g)$ is
  either nearly parallel or conformally equivalent to the standard
  7-sphere.
\end{theorem}
\begin{theorem}\label{main2}
  Let $(M,g,J)$ be a compact 6-dimensional manifold locally
  conformally equivalent to a nearly K\"ahler manifold.  Then
  $(M,g,J)$ has constant scalar curvature if and only if $(M,g)$ is
  either nearly parallel or conformally equivalent to the 6-sphere
  with standard metric.
\end{theorem}

In the last section we characterize complete Einstein $G_2$ and
$SU(3)$ manifolds in the strict class $\mathcal X_1+\mathcal X_4$ and
$\mathcal W_1+\mathcal W_4$, respectively.  The phrases `strict class'
is used here to indicate that the $G$-structure is not in any
sub-class of the one given. So a $G_2$-structure strictly in class
$\mathcal X_1$ must, in particular, have non-trivial intrinsic
torsion.

The results obtained in this paper are direct consequences of the
following.  The $G$-structures under consideration are described by
the existence of certain \emph{fundamental differential forms} \(
\om1,\dots,\om{p} \), whose exterior derivatives determine the
corresponding intrinsic torsion in full.  First order identities on
the $G$-invariant components of the intrinsic torsion descend from the
closure of $d\omega_i$.  These equations in general have non-trivial
consequences as is seen by the examples considered here.

The relations coming from the second derivatives of the forms may also
be seen as consequences of the first Bianchi identity, see for
instance~\cite{S&C,math.DG/0305124}.

\begin{acknowledgements}
  The authors wish to thank Andrew Swann, Ilka Agricola and Thomas
  Friedrich for helpful conversations and Simon Chiossi for many
  helpful comments.  The first named author was supported by the
  Junior Research Group ``Special Geometries in Mathematical
  Physics'', of the Volkswagen foundation and the SFB~647
  ``Space--Time--Matter'', of the DFG.  The authors wish to thank UC
  Riverside for hospitality and Yat Sun Poon for support during the
  initial stage of this project.
\end{acknowledgements}

\section{A lemma}
The key to obtaining the results is the observation
\begin{lemma}\label{lem:1}
  Let $M$ be a connected manifold equipped with a differentiable
  function $\phi\not\equiv 0$ and a one-form $\alpha$ such that
  \begin{gather}
    \label{eq:1}
    \begin{split}
    d\alpha &= 0,\\
    d\phi+\phi&\alpha = 0.
    \end{split}
  \end{gather}
  Then $\phi$ is nowhere zero and $\alpha = -d\log\abs{\phi}$.
\end{lemma}
\begin{proof}
  Let $\phi$ and $\alpha$ be a function and one-form as in
  equation~\eqref{eq:1}.  By Poincar\'e's Lemma we may choose a
  covering $U_i$ of $M$ and functions $f_i\colon U_i\to\mathbb R$ such
  that $\alpha\lvert_{U_i} = df_i$.  Then equation~\eqref{eq:1}
  implies that the product $\phi\exp(f_i)$ is constant over each
  $U_i$.  Therefore, if $\phi(p)\not= 0$ at some point $p$ in, say
  $U_0$, then $\phi \not= 0$ over all $U_0$ and therefore also on each
  $U_j$ that overlaps $U_0$.  The conclusion now follows from
  connectedness.
\end{proof}

\section{The $G_2$ case}
\label{sec:g_2-case}

A $G_2$-manifold is a $7$-dimensional manifold $M$ equipped
with a special, so-called fundamental three-form $\phi$, required to
satisfy the following non-degeneracy condition
\begin{equation}
  i_X\phi\bw i_Y\phi\bw\phi=6g(X,Y)\vol(g),\label{eq:21}
\end{equation}
for some Riemannian metric $g$ and orientation on $M$.  The notation
$i_X\phi$ means interior product of the vector field $X$ with the
three-form $\phi$.  It is well known that the covariant derivative of
the fundamental three-form is determined by the exterior derivatives
of $\phi$ and its Hodge dual $\Hodge\phi$.  Using the representation
theory of $G_2$ on the exterior algebra one may write these
differentials as
\begin{gather*}
  d\phi=\tau_0\Hodge\phi+3\tau_1\bw\phi+\Hodge\tau_3,\\
  d\Hodge\phi=4\tau_1\bw\Hodge\phi+\tau_2\bw\phi,
\end{gather*}
for suitable forms $\tau_p\in\Omega^p$.  In terms of the $G_2$
invariant splittings of the exterior algebra, $\tau_0 \in
\Omega_{1}^0, ~\tau_1 \in \Omega^1_7, ~\tau_2 \in \Omega_{14}^2,
~\tau_3 \in \Omega^3_{27}$.  The notation $\Omega^p_d$ indicates the
space of $p$-forms taking values in the $d$-dimensional $G_2$
irreducible subspace $\Lambda^p_d\subset\Lambda^p$. The one-form
$\tau_1$ is also known as the \emph{Lee form} of the $G_2$ manifold.
The forms $\tau_0,\tau_1,\tau_2,\tau_3$ correspond to the
Fern\'andez-Gray classes~\cite{MR696037}as follows
\begin{gather*}
  \tau_0\leftrightarrow \mathcal X_1,\quad
  \tau_2\leftrightarrow \mathcal X_2,\quad
  \tau_3\leftrightarrow \mathcal X_3,\quad
  \tau_1\leftrightarrow \mathcal X_4.
\end{gather*}
When we speak of the intrinsic torsion $\tau$ of a $G_2$ structure we
mean the form of mixed degree $\tau=\tau_0+\tau_1+\tau_2+\tau_3$ fixed
by the fundamental three-form as above.

In particular, $G_2$-manifolds in the class $\mathcal X_1$ are
characterized by the conditions $\tau_1=\tau_2=\tau_3=0$ and are
called nearly parallel $G_2$-manifolds. It is well known that these
spaces are Einstein with positive scalar curvature.  From this it
follows that $\tau_0$ is constant~\cite{BFKM,GSal}.

A $G_2$-manifold in the Fern\'andez-Gray class $\mathcal X_1+\mathcal
X_4$ satisfies $\tau_2=\tau_3=0$.  The structure equations for this case
reduce to
\begin{gather*}
  d\phi=\tau_0\Hodge\phi+3\tau_1\bw\phi,\\
  d\Hodge\phi=4\tau_1\bw\Hodge\phi,
\end{gather*}
from which one infers
\begin{equation}
  \label{eq:18}
  \begin{array}{c}
  d^2\phi=(d\tau_0+\tau_0\tau_1)\bw\Hodge\phi+3d\tau_1\bw\phi,\\
  d^2\Hodge\phi=4d\tau_1\bw\Hodge\phi.
  \end{array}
\end{equation}
The latter equation implies that the component $(d\tau_1)_{7} \in
\Omega^2_{7}$ vanishes.  Using this in the first equation, one deduces
that the complementary component $(d\tau_1)_{14} \in \Omega^2_{14}$
also vanishes.  Thus we recover the fact (see~\cite{MR1386249}) that
$G_2$ structures in this class are, locally, conformal to a nearly
parallel structure. Observe that equations~\eqref{eq:18} furthermore
give us $d\tau_0+\tau_0\tau_1=0$.  Now our Lemma~\ref{lem:1} applies
with $\phi=\tau_0$ and $\alpha=\tau_1$.  The connectedness of $M$
leads to the conclusion in the form of this
\begin{theorem}\label{thm:g2}
  Suppose $M$ is a $7$-dimensional manifold with a $G_2$-structure
  $\phi$ in the Fern\'andez-Gray class $\mathcal X_1 + \mathcal X_4$.
  Then $\phi$ is either of class $\mathcal X_4$, in which case
  $(M,\phi)$ is locally conformal to a parallel $G_2$-manifold, or
  $(M,\phi)$ is conformal to a nearly parallel $G_2$-manifold.
\end{theorem}

\section{The $6$-dimensional almost Hermitian case}

An almost Hermitian manifold is a Riemannian manifold $(M^{2m},g)$
equipped with an orthogonal almost complex structure $J$.  The metric
and the almost complex structure then define the fundamental two-form
of the almost Hermitian structure:
\begin{gather*}
  \omega(X,Y):=g(JX,Y).
\end{gather*}
As opposed to the $G_2$ case above and the case of $\SU(3)$ below, the
components of the intrinsic torsion of a Hermitian structure cannot
all be identified with differential forms~\cite{MR1004008}.  Instead,
the intrinsic torsion is detected by $d\omega$ along with the
Nijenhuis tensor $N_J$.  

From now on $m$ is taken to be at least $3$.

\subsection{The Nijenhuis Tensor}
\label{sec:3.6}

For an almost complex structure $J$ the Nijenhuis tensor measures the
failure of the eigenspaces $J$ in the complexified tangent space to be
involutive.  Let $\pi'(X)=\frac12(X-iJX)$ and $\pi''=\frac12(X+iJX)$
be the projections to the $i$-eigenspace $T'$ and the $-i$-eigenspace
$T''$, respectively.  We set
\begin{align*}
  N_J(X,Y):=&\pi'[\pi''X,\pi''Y]+\pi''[\pi'X,\pi'Y]\\
  =& \frac14\left([X,Y]-[JX,JY]+J[JX,Y]+J[X,JY]\right).
\end{align*}
Using the metric we obtain an algebraically equivalent $3$-tensor
\begin{equation*}
  N_J(X,Y;Z)=g(N_J(X,Y),Z),
\end{equation*}
with the property $N_J(JX,Y;Z)=N_J(X,JY;Z)=N_J(X,Y;JZ)$.
Equivalently,
\begin{equation*}
  N_J\in\real{\Lambda^{(2,0)}\otimes\Lambda^{(1,0)}}.
\end{equation*}
See for instance~\cite{MR1299398} for an explanation of the notation.

The space $\real{\Lambda^{(3,0)}}$ is a subspace of
$\real{\Lambda^{(2,0)} \otimes \Lambda^{(1,0)}}$ in the natural way.
The projection $\real{\Lambda^{(2,0)} \otimes \rep{(1,0)}{}} \to
\real{\Lambda^{(3,0)}}$ is given simply by skew-symmetrization.  Write
$V$ for the orthogonal complement of $\real{\Lambda^{(3,0)}}$ in
$\real{\Lambda^{(2,0)} \otimes \Lambda^{(1,0)}}$.  Then we may split
the Nijenhuis tensor accordingly
\begin{equation*}
  N_J=N_J^{3,0}+N_J^V.
\end{equation*}
One may now deduce that
\begin{equation}\label{eq:23}
  (d\omega)^{3,0}(X,Y,Z)
  =3g(N_J^{3,0}(X,Y),JZ)=3N_J^{3,0}(X,Y,JZ).
\end{equation}
The structure equations for an almost Hermitian manifold now can be
written
\begin{equation}
  \label{eq:19}
  \begin{array}{c}
    d\omega=-3J_{(1)}N_J^{3,0}+2\sigma_1\bw\omega+\sigma_3,\qquad
    N_J=N_J^{3,0}+N^V_J.
  \end{array}
\end{equation}
The first equation here employs the conventions
$(J_{(1)}\alpha)(X,Y,\dots) := -\alpha(JX,Y,\dots)$,
see~\cite{MR2162420}.  This action of the complex structure $J$ on
differential forms is, generally speaking, distinct from the usual
action given by $(J\alpha)(X_1,\dots,X_p) :=
(-1)^p\alpha(JX_1,\dots,JX_p)$.

The Gray-Hervella classes of an almost Hermitian manifold are in the
following correspondence with the components in~\eqref{eq:19}
\begin{gather*}
  N_J^{3,0}\leftrightarrow \mathcal W_1,\quad
  N^V_J\leftrightarrow \mathcal W_2,\quad
  \sigma_3\leftrightarrow \mathcal W_3,\quad
  \sigma_1\leftrightarrow \mathcal W_4.
\end{gather*}
Almost Hermitian manifolds in the class $\mathcal W_1$, called nearly
K\"ahler manifolds, are characterized by the conditions
$N_J^V=\sigma_1=\sigma_3=0$ or equivalently, by demanding that the
covariant derivative of the almost complex structure with respect to
the Levi-Civita connection be skew-symmetric, $(\LC_XJ)X=0$ \cite{Gr}.
\begin{remark}
  Almost Hermitian manifolds in the class $\mathcal W_1+\mathcal
  W_3+\mathcal W_4$ are characterized by $N_J^V=0$, i.e. the Nijenhuis
  tensor is totally skew-symmetric. This amounts to the existence of a
  linear connection preserving the almost Hermitian structure with
  totally skew-symmetric torsion \cite{FI}. In dimension $6$, this
  class can also be characterized by the property that the Nijenhuis
  tensor is either everywhere non-degenerate (of constant signature)
  or vanishes identically \cite{Br2}. These manifolds are called
  quasi-integrable and investigated in details in \cite{Br2}.
\end{remark}
For an arbitrary one-form the relation
\begin{equation*}
  (d\theta)(X,Y) - (d\theta)(JX,JY) -
    d(J\theta)(JX,Y) + Jd(J\theta)(JX,Y) 
    = 4g(N_J(X,Y),J(J\theta)^\#)
\end{equation*}
holds.  Writing $d\theta^{2,0}=\frac12(d\theta-Jd\theta)$ for the
projection of $d\theta$ to $\real{\Lambda^{(2,0)}}$ we have
\begin{lemma}
  \label{lem:2} Suppose $(g,J)$ is an almost Hermitian structure
  in class $\mathcal W_1+\mathcal W_3+\mathcal W_4$.  Let $\theta$ be
  a one-form and write $\theta':=J\theta$.  Then
  \begin{align*}
    (d\theta)^{2,0} + J_{(1)}(d\theta')^{2,0} &=
    \tfrac23\theta'\hook(d\omega)^{3,0}
  \end{align*}
\end{lemma}

\subsection{$\mathbf{SU(3)}$-structures}
\label{sec:spec-almost-herm}
A $6$-dimensional manifold with an $\SU(3)$-structure comes equipped
with data $(g,J,\omega,\psi_+,\psi_-)$ invariant with respect to the
action of $\SU(3)$. Here $g$ is a Riemannian metric, $J$ is an almost
complex structure, $\omega$ the fundamental two-form and $\psi_+$ and
$\psi_-$ are three-forms such that $\Psi:=\psi_++i\psi_-$ is a complex
$(3,0)$-form.  These invariant tensors are not independent, in fact
the triple $(\omega,\psi_+,\psi_-)$ with $\psi_++i\psi_-$ decomposable
and compatible with $\omega$ by means of the equations below, defines
both $g$ and $J$, see~\cite{hitchin97:_lagran}.  Clearly the triple
$(g,\psi_+,\psi_-)$ will do the same.  We choose a normalization with
the following relations
\begin{equation}
  \begin{array}{c}
    \omega(X,JY)=g(X,Y),\\
    \omega\bw\psi_+=0=\omega\bw\psi_-,\\
    3\psi_+\bw\psi_-=2\omega^3=12\vol_g,\\
    \Hodge\omega=\tfrac12\omega^2,\quad
    \Hodge\psi_+=\psi_-,\quad
    J\psi_+=-\psi_-.
  \end{array}
\end{equation}

\subsubsection{Torsion classes and structure equations}
\label{sec:2}

Under the action of $\SU(3)$, $\Lambda^{(3,0)}=\mathbb C$ and
$V\otimes\mathbb C\cong\Lambda^{(1,1)}_0$. This means that $V\cong
2\su(3)$ and $\real{\Lambda^{(3,0)}}\cong2\mathbb R$.  Moreover, for
an $\SU(3)$-structure $(\omega ,\psi_\pm)$, the components of the
Nijenhuis tensor can be computed from components of
$(d\omega,d\psi_\pm)$.  In fact there are algebraic correspondences
\begin{gather*}
  N_J^{(3,0)} \leftrightarrow (d\omega)^{(3,0)+(0,3)} \leftrightarrow
  \left((d\psi_+)^{(0,0)},(d\psi_-)^{(0,0)}\right),\\
  N_J^V \leftrightarrow \left((d\psi_+)^{(2,2)}_0,(d\psi_-)^{(2,2)}_0\right).
\end{gather*}
The first arrow is given by equation~\eqref{eq:23}, the second is
described below. The remaining one has a similar description which we
will not need here.

The $\SU(3)$-structure function $\nabla^g\Psi$ is completely
determined by the exterior derivatives of the three forms
$\omega,\psi_+$ and $\psi_-$.  These may be written as
\begin{gather}
  d\omega = 3\left(\sigma_0^+\psi_+ - \sigma_0^- \psi_-\right) +
  2\sigma_1^+\bw\omega + \sigma_3,\notag\\
  \label{eq:9}
  d\psi_+ = - 2\sigma_0^- \omega^2 + 3\sigma_1^+\bw\psi_+ -
  \sigma_1^-\bw\psi_- +
  \sigma_2^+\bw\omega,\\
  d\psi_- = - 2\sigma_0^+\omega^2 + 3\sigma_1^+\bw\psi_- +
  \sigma_1^-\bw\psi_+ + \sigma_2^-\bw\omega.\notag
\end{gather}
Here $\sigma_p^{\pm}$ are $p$-forms and $\sigma_3$ is a three-form.
They correspond roughly to the classes $\mathcal W_1^+,~\mathcal
W_1^-,~\mathcal W_4,~\mathcal W_5,~\mathcal W_2^+,~\mathcal W_2^-,$
and $\mathcal W_3$ of~\cite{MR1922042}, respectively (see
also~\cite{bor01:_bochn_g}).  These determine the Gray-Hervella
classes of the underlying almost Hermitian structure in the obvious
way.

\begin{remark}
  The one-form $\sigma_1^+$ is, in fact, the Lee form of the almost
  Hermitian structure. In contrast $\sigma_1^-$ is conformally
  invariant.  Therefore $\sigma_1^-$ does not really correspond to the
  class $\mathcal W_5$ but rather to ``$3W_4+2W_5$''.  This choice for
  the one-forms was introduced by Mart\'in Cabrera~\cite{MR2162420}.
\end{remark}

\subsubsection{A transformation}
\label{sec:3}

Set $\lambda:=\sigma_0^++i\sigma_0^- $ and $\Lambda:=\abs{\lambda}$.
In neighbourhoods with $\lambda$ non-vanishing an argument
$\varphi:=\arg(\lambda):=\arctan\left(\frac{\sigma_0^-
  }{\sigma_0^+}\right)$ may be chosen.  We then set
\begin{gather*}
  \tilde\omega:=\Lambda^2\omega,\\
  \tilde\psi_+:=\Lambda^2\left(\sigma_0^+\psi_+ - \sigma_0^- \psi_-\right),\\
  \tilde\psi_-:=\Lambda^2\left(\sigma_0^-\psi_+ + \sigma_0^+ \psi_-\right).
\end{gather*}
This gives the somewhat simpler structure equations
\begin{gather*}
  d\tilde\omega := 3\tilde\psi_+ + 2{\tilde\sigma_1^+}\bw\omega + \tilde\sigma_3,\\
  d\tilde\psi_+ := 3{\tilde\sigma_1^+}\bw\tilde\psi_+ -
  {\tilde\sigma_1^-}\bw\tilde\psi_- + {\tilde\sigma_2^+}\bw\tilde\omega{},\\
  d\tilde\psi_- := - 2{\tilde\omega}^2 +
  3{\tilde\sigma_1^+}\bw\tilde\psi_- +
  {\tilde\sigma_1^-}\bw\tilde\psi_- +
  {\tilde\sigma_2^-}\bw\tilde\omega{}.
\end{gather*}
where
\begin{gather*}
  {\tilde\sigma_1^+}:=\sigma_1^+ +\Lambda^{-1}d\Lambda,\qquad
  {\tilde\sigma_1^-} :=
  \sigma_1^- - d\varphi,\\
  {\tilde\sigma_2^+}:={\sigma_0^+}\sigma_2^+-{\sigma_0^-}\sigma_2^-,
  \qquad {\tilde\sigma_2^-}:=
  {\sigma_0^-}\sigma_2^++{\sigma_0^+}\sigma_2^-,\\
  \tilde\sigma_3:=\Lambda^2\sigma_3.
\end{gather*}
In particular, the structure equations of a nearly K\"ahler
$6$-manifold can always be put on the form~\cite{Hitchin, Sal}
$$d\omega = 3\psi_+, \qquad d\psi_- = - 2\omega^2.$$
It is well known that these spaces are Einstein with positive
scalar curvature \cite{MR0293537}.

\section{Locally conformally nearly K\"ahler $6$-folds}
\label{sec:4}

For a $6$-dimensional almost Hermitian manifold in the class $\mathcal
W_1+\mathcal W_4$ a (possibly local) choice of trivialization
$(\psi_+,\psi_-)$ allows us to write the structure
equations~\eqref{eq:19}, \eqref{eq:9} as
\begin{gather}
  d\omega = 3\left(\sigma_0^+\psi_+ - \sigma_0^- \psi_-\right) +
  2\sigma_1^+\bw\omega,\notag\\
  \label{eq:10}
  d\psi_+ = - 2\sigma_0^- \omega^2 + 3\sigma_1^+\bw\psi_+ -
  \sigma_1^-\bw\psi_-,\\
  d\psi_- = - 2\sigma_0^+\omega^2 + 3\sigma_1^+\bw\psi_- +
  \sigma_1^-\bw\psi_+.\notag
\end{gather}
Differentiating each equation yields
\begin{gather}
  \label{eq:11}
  0 = 3(d\sigma_0^+ + \sigma_0^+\sigma_1^+ - \sigma_0^-  \sigma_1^-)\psi_+ -
  3(d\sigma_0^-  + \sigma_0^- \sigma_1^+ + \sigma_0^+\sigma_1^-)\psi_- +
  2d\sigma_1^+\omega,\\
  \label{eq:12}
  0 = - 2(d\sigma_0^- + \sigma_0^- \sigma_1^+ +
  \sigma_0^+\sigma_1^-)\omega^2 + 3d\sigma_1^+\psi_+ -
  d\sigma_1^-\psi_- ,\\ 
  \label{eq:13}
  0 = - 2(d\sigma_0^+ + \sigma_0^+\sigma_1^+ - \sigma_0^-
  \sigma_1^-)\omega^2 + 3d\sigma_1^+\psi_- + d\sigma_1^-\psi_+. 
\end{gather}
These have the following immediate consequences.
Equation~\eqref{eq:11} shows that $d\sigma_1^+$ is a $(2,0)+(0,2)$
form as a linear combination of one-forms contracted with $\psi_+$ and
$\psi_-$.  Using standard identities such as $J(\sigma\bw\psi_+) =
\sigma\bw\psi_-$ for an arbitrary two-form $\sigma$ and
$J(\sigma\bw\omega)=\sigma\bw\omega$ for a $(1,1)$-form, as well as
$\Hodge(\theta\bw\psi_-) = \theta\hook\psi_+ = (J\theta)\hook\psi_-$
for a one-form $\theta$, leads to the equivalent set of equations:
\begin{gather}
  \label{eq:14}
  d\sigma_0^+ + \sigma_0^+\sigma_1^+ - \sigma_0^- \sigma_1^- = J(d\sigma_0^-
  + \sigma_0^- \sigma_1^+ + \sigma_0^+\sigma_1^-),\\
  \label{eq:15}
  d\sigma_1^+ = 3(d\sigma_0^- + \sigma_0^- \sigma_1^+ +
  \sigma_0^+\sigma_1^-)\hook\psi_+,\\
  \label{eq:16}
  (d\sigma_1^-)^{2,0}=7(d\sigma_0^- + \sigma_0^- \sigma_1^+ +
  \sigma_0^+\sigma_1^-) \hook\psi_-.
\end{gather}

\begin{lemma}\label{lem:l.c}
  Suppose $(M^6,\omega,J)$ is an almost Hermitian manifold in the
  class $\mathcal{W}_1+\mathcal W_4$.  Then the Lee form is closed if
  and only if $(\omega,J)$ is either \emph{globally} conformal to a
  nearly K\"ahler structure $(\omega',J')$ on $M$ or locally
  conformally equivalent to a K\"ahler structure.
\end{lemma}
\begin{proof}
  Suppose the Lee-form is closed $d\sigma_1^+=0$.  Locally, we pick a
  smooth trivialisation $(\psi_+,\psi_-)$ of $\real{\Lambda^{(3,0)}}$.
  Then, locally, equations~\eqref{eq:15} and \eqref{eq:14} show that
  \begin{gather*}
    d\sigma_0^+ + \sigma_0^+\sigma_1^+ - \sigma_0^- \sigma_1^- = 0,\\
    d\sigma_0^-  + \sigma_0^- \sigma_1^+ + \sigma_0^+\sigma_1^- = 0,
  \end{gather*}
  whence
  \begin{gather*}
    d\left((\sigma_0^+)^2 + (\sigma_0^-)^2\right)+((\sigma_0^+)^2 +
    (\sigma_0^-)^2)(2\sigma_1^+)=0.
  \end{gather*}
  However,
  \begin{equation*}
    \phi := (\sigma_0^+)^2 + (\sigma_0^-)^2 = \tfrac19
    \norm{d\omega^{3,0}}^2
  \end{equation*}
  is a globally well-defined, smooth function, and $\alpha:=2\sigma_1^+$
  is closed.  So Lemma~\ref{lem:1} applies and we conclude that
  $\norm{d\omega^{3,0}}$ is either non-zero everywhere, or it vanishes
  at all points.
\end{proof}
\begin{remark}
  Note that in dimension $6$ the Nijenhuis tensor is either everywhere
  non-degenerate (of constant signature) or vanishes identically not
  only for $\mathcal W_1+\mathcal W_4$, but for the whole class
  $\mathcal W_1+\mathcal W_3+\mathcal W_4$ \cite{Br2}.
\end{remark}

\begin{theorem}\label{thnkel}
  Let $M$ be a $6$-dimensional manifold with an almost Hermitian
  structure $(\omega,J)$ in the Gray-Hervella class $\mathcal
  W_1+\mathcal W_4$.  Then either $(\omega,J)$ is locally conformally
  equivalent to K\"ahler structure on $M$, or $(\omega,J)$ is a
  conformal transformation of a nearly K\"ahler structure.
\end{theorem}
\begin{proof}
  Write $M$ as a disjoint union $M_0\cup M_1$ where
  \begin{gather*}
    M_0:=\set{x\in M}{(d\omega)^{3,0}=0},\qquad
    M_1:=\set{x\in M}{(d\omega)^{3,0}\not=0}.
  \end{gather*}
  On the open submanifold $M_1$ there is a canonical choice of
  trivialization of $\real{\Lambda^{(3,0)}}$ given by taking $\psi_+=
  (d\omega)^{3,0}$.  After a suitable transformation (as in
  section~\ref{sec:3}) we obtain the structure equations
  \begin{gather}
    d\tilde\omega = 3\tilde\psi_+ +
    2\tilde\sigma_1^+\bw\tilde\omega,\notag\\
    d\tilde\psi_+ = 3\tilde\sigma_1^+\bw\tilde\psi_+ -
    \tilde\sigma_1^-\bw\tilde\psi_-,\\
    d\tilde\psi_- = - 2\tilde\omega^2 +
    3\tilde\sigma_1^+\bw\tilde\psi_- +
    \tilde\sigma_1^-\bw\tilde\psi_+.\notag 
  \end{gather}
  Equations~\eqref{eq:14}, \eqref{eq:15} and~\eqref{eq:16} then become
  \begin{gather}
    \tilde\sigma_1^+ = J\tilde\sigma_1^-,\notag\\
    d\tilde\sigma_1^+ = 3\tilde\sigma_1^-\hook\tilde\psi_+,\label{eq:2}\\
    (d\tilde\sigma_1^-)^{2,0}=7\tilde\sigma_1^-\hook\tilde\psi_-\notag.
  \end{gather}
  Using Lemma~\ref{lem:2} with
  $\theta=\tilde\sigma_1^+,~\theta'=\tilde\sigma_1^-$, and the
  identity $J_{(1)}(\sigma\hook\psi_{\pm}) = (J\sigma)\hook\psi_{\pm}
  = \mp\sigma\hook\psi_{\mp}$ valid for all one-forms $\sigma$, we get
  \begin{equation*}
    d\tilde\sigma_1^+ - J_{(1)}(d\tilde\sigma_1^-)^{2,0} =
    -2\tilde\sigma_1^+\hook\tilde\psi_- =
    -2\tilde\sigma_1^-\hook\tilde\psi_+.
  \end{equation*}
  This is only compatible with the relations~\eqref{eq:2} if
  $\tilde\sigma_1^-\hook\tilde\psi_+=0$.  Therefore $\tilde\sigma_1^+
  =\tilde\sigma_1^-=0$ and the original one-forms $\sigma_1^\pm$ are,
  in fact, exact on $M_1$.  Moreover, on the interior $M^o_0$ of
  $M_0$, $d\omega = 2\sigma_1^+\bw\omega$, so
  $d\sigma_1^+\lvert_{M_0^o}=0$ also holds.
  
  So the set of points at which $d\sigma_1^+\not=0$, which clearly is
  open, is the common boundary of two open sets in $M$, at least one
  of which is non-empty.  Therefore $d\sigma_1^+=0$ on all of $M$ and
  Lemma~\ref{lem:l.c} completes the proof.
\end{proof}

\section{Proof of Theorem~\ref{main1} and Theorem~\ref{main2}}

Theorem~\ref{thm:g2} and Theorem~\ref{thnkel} show that the Riemannian
manifold $(M,g)$ is globally conformal to an Einstein space of
positive scalar curvature. Further, if the scalar curvature is
constant then the Obata Theorem (see the proof in \cite{LP}) tells us
that the conformal transformation making the metric Einstein is
trivial, or else $(M,g)$ is the standard sphere. \qed

\section{Curvature classification}
The Riemannian curvature tensor of a nearly K\"ahler $6$-manifolds or
a nearly parallel $G_2$-manifold is especially simple.  In fact,
viewing curvature tensors as bundle endomorphisms $R \colon \Lambda^2
\to \Lambda^2$, the curvature splits as
\begin{equation}
  \label{eq:22}
  \RC=R^{\lie g}+\frac{s_g}{2n(n-1)}\Id_{\Lambda^2}.
\end{equation}
where $R^{\lie g}$, formally, is the curvature tensor of a space with
holonomy algebra $\lie g$ and $n$ is the dimension of the underlying
space.  This formula is reminiscent of the curvature formula for a
Riemannian manifold with holonomy $\SP(n)\SP(1)$ of~\cite{MR1004008}.
In the cases of concern, $\lie g$ and $n$ are equal to $\su(3)$ and
$6$ for nearly K\"ahler and $\lie g_2$ and $7$ for nearly parallel
$G_2$.  In either situation the tensor $R^{\lie g}$ takes values in a
$G$-irreducible subspace of the space of algebraic Weyl tensors, i.e.,
algebraic curvature tensors with vanishing Ricci contraction.
Standard identities~\cite{Besse} now make it possible to deduce the
form of the curvature tensor for an almost Hermitian or $G_2$ space of
type $\mathcal W_1+\mathcal W_4$, or $\mathcal X_1+\mathcal X_4$,
respectively.  The conformal invariance of the Weyl curvature implies
that a manifold in the more general classes are Einstein if and only
if their curvature is of the same form of nearly parallel structures.

\begin{theorem} 
  \begin{enumerate}[(a)]
  \item Suppose $(M,\phi)$ is a $G_2$ manifold of strict type
    $\mathcal X_1+\mathcal X_4$ such that the associated metric $g$ is
    complete and Einstein.  Then $(M,g)$ is isometric to either the
    sphere, the hyperbolic space or the euclidean space equipped with
    a constant curvature metric.  Furthermore, the $G_2$ structure
    $\phi$ is a conformal change of a suitable restriction of the
    standard nearly parallel structure on the $7$ dimensional sphere.
  \item Suppose $(M,\omega,J)$ is an almost Hermitian $6$-manifold of
    strict type $\mathcal W_1+\mathcal W_4$ such that the associated
    metric $g$ is complete and Einstein.  Then $(M,g)$ is isometric to
    either the sphere, hyperbolic space or euclidean space equipped
    with a constant curvature metric.  Furthermore, the almost
    Hermitian structure $(\omega,J)$ is a conformal change of a
    suitable restriction of the standard nearly K\"ahler structure on
    the $6$ dimensional sphere.
  \end{enumerate}
\end{theorem}
\begin{proof}
  This is an immediate consequence of Theorem~\ref{thm:g2},
  Theorem~\ref{thnkel}, uniqueness of the nearly parallel structure
  compatible with the round metric on the $7$-sphere, uniqueness of
  the nearly K\"ahler structure compatible with the round metric on
  the $6$-sphere, see Friedrich~\cite{Friedrich} and the Main Theorem
  of~\cite{Kuhnel}.
\end{proof}
\begin{corollary}
  Suppose $M$ is a $G_2$ manifold or almost Hermitian 6-manifold of
  strict type $\mathcal X_1+\mathcal X_4$ or $\mathcal W_1+\mathcal
  W_4$, respectively. Assume that the Riemannian curvature of $M$ is
  of the form \eqref{eq:22}.  Then $M$ has constant sectional
  curvature and in particular $R^{\lie g}=0$.
\end{corollary}

\bibliographystyle{plain}

\end{document}